\documentclass[smallextended,envcountsect]{svjour3}
\smartqed
\usepackage{lineno}
\modulolinenumbers[5]

\journalname{Journal of XXXX}
\usepackage{amssymb,amsmath,amsfonts}
\usepackage{epsfig}
\usepackage[bookmarksnumbered, colorlinks, urlcolor=blue, linkcolor=blue, citecolor=red, plainpages=false, pdfstartview=FitW]{hyperref}
\usepackage{mathptmx}
\usepackage{diagbox}
\usepackage{graphicx}
\usepackage{color}
\usepackage{caption}
\usepackage{wrapfig}
\usepackage{graphics}
\usepackage{enumerate}
\usepackage{amsxtra}
\usepackage{latexsym}
\usepackage{multirow}
\usepackage{subfigure}
\usepackage{epstopdf}
\usepackage{pdfpages}
\usepackage{algorithm}
\usepackage{algorithmic}
\usepackage{verbatim}
\usepackage{cite}
\usepackage{url}
\newcommand{\thmlist}{
\begin{list}{Step 1}
{\setlength{\leftmargin}{0.6 in}\setlength{\labelwidth} {0.5 in}}}

\newcommand{\alglist}{
\begin{list}{Step 1}
{\setlength{\leftmargin}{1.1 in} \setlength{\labelwidth}{1.0 in}}}

 \renewcommand{\proof} {\noindent {\bf Proof.} \quad}

\renewcommand{\subtitle}[1]{\color{blue}}

\def\red#1{\color{red}{#1}\color{black}}
\def\blue#1{\color{blue}{#1}\color{black}}

\begin{document}

\title{Primal-dual path-following methods and the trust-region updating
strategy for linear programming with noisy data}
\titlerunning{Primal-dual path-following methods and the trust-region updating strategy}
\author{Xin-long Luo \textsuperscript{$\ast$} \and Yi-yan Yao}
\authorrunning{Luo \and Yao}

\institute{
   Xin-long Luo, Corresponding author
   \at
   School of Artificial Intelligence,
   Beijing University of Posts and Telecommunications, P. O. Box 101,
   Xitucheng Road  No. 10, Haidian District, 100876, Beijing China,
   \email{luoxinlong@bupt.edu.cn}
   \and
   Yi-yan Yao
   \at
   School of Artificial Intelligence,
   Beijing University of Posts and Telecommunications, P. O. Box 101,
   Xitucheng Road  No. 10, Haidian District, 100876, Beijing China,
   \email{yaoyiyan@bupt.edu.cn}
}

\date{Received: date / Accepted: date}
\maketitle

\begin{abstract}
  In this article, we consider the primal-dual path-following method and
  the trust-region updating strategy for the standard linear programming problem.
  For the rank-deficient problem with the small noisy data, we also give
  the preprocessing method based on the QR decomposition with column pivoting.
  Then, we prove the global convergence of the new method when the initial
  point is strictly primal-dual feasible. Finally, for some rank-deficient
  problems with or without the small noisy data from the NETLIB collection,
  we compare it with other two popular interior-point methods, i.e. the
  subroutine pathfollow.m and the built-in subroutine linprog.m of the
  MATLAB environment. Numerical results show that the new  method is more
  robust than the other two methods for the rank-deficient problem with
  the small noise data.
\end{abstract}


 \vskip 2mm



\section{Introduction}

\vskip 2mm

In this article, we are mainly concerned with the linear programming problem
with the small noisy data as follows:
\begin{align}
  & \min_{x \in \Re^{n}} c^{T}x, \hskip 2mm \text{subject to} \; Ax = b,
  \; x \ge 0,     \label{LPND}
\end{align}
where $c$ and $x$ are vectors in $\Re^{n}$, $b$ is a vector in $\Re^{m}$,
and $A$ is an $m \times n$ matrix. For the problem \eqref{LPND}, there are many
efficient methods to solve it such as the simplex methods \cite{Pan2010,ZYL2013},
the interior-point methods \cite{FMW2007,Gonzaga1992,NW1999,Wright1997,YTM1994,Zhang1998}
and the continuous methods \cite{AM1991,CL2011,Monteiro1991,Liao2014}. Those
methods are all assumed that the constraints of problem \eqref{LPND} are consistent,
i.e. $rank(A,\, b) = rank(A)$. For the consistent system of redundant constraints,
references \cite{AA1995,Andersen1995,MS2003} provided a few preprocessing
strategies which are widely used in both academic and commercial linear programming
solvers.

\vskip 2mm

However, for a real-world problem, since it may include the redundant constraints
and the measurement errors, the rank of matrix $A$ may be deficient and the
right-hand-side vector $b$ has small noise. Consequently, they may lead to
the inconsistent system of constraints \cite{AZ2008,CI2012,LLS2020}.
On the other hand, the constraints of the original real-world problem are
intrinsically consistent. Therefore, we consider the least-squares
approximation of the inconsistent constraints in the linear programming problem
based on the QR decomposition with column pivoting. Then, according to the
first-order KKT conditions of the linear programming problem, we convert the
processed problems into the equivalent problem of nonlinear equations with
nonnegative constraints. Based on the system of nonlinear equations with
nonnegative constraints, we consider a special continuous Newton flow with
nonnegative constraints, which has the nonnegative steady-state solution
for any nonnegative initial point. Finally, we consider a primal-dual
path-following method and the adaptive trust-region updating strategy to follow
the trajectory of the continuous Newton flow. Thus, we obtain an optimal
solution of the original linear programming problem.

\vskip 2mm

The rest of this article is organized as follows. In the next section, we
consider the primal-dual path-following method and the adaptive trust-region
updating strategy for the linear programming problem. In section 3, we analyze
the global convergence of the new method when the initial point is strictly
primal-dual feasible. In section 4, for the rank-deficient problems with or
without the small noise, we compare the new method with two other popular
interior-point methods, i.e. the traditional path-following method
(pathfollow.m in p. 210, \cite{FMW2007}) and the predictor-corrector algorithm
 (the built-in subroutine linprog.m of the MATLAB environment
 \cite{MATLAB,Mehrotra1992,Zhang1998}). Numerical results show that
the new method is more robust than the other two methods for the rank-deficient
problem with the small noisy data. Finally, some discussions are given
in section 5. $\|\cdot\|$ denotes the Euclidean vector norm or its induced
matrix norm throughout the paper.

\vskip 2mm

\section{Primal-dual path-following methods and the trust-region updating strategy}

\subsection{The continuous Newton flow} \label{SUBSECNF}

For the linear programming problem \eqref{LPND}, it is well known that its optimal
solution $x^{\ast}$ if and only if it satisfies the following Karush-Kuhn-Tucker
conditions (pp. 396-397, \cite{NW1999}):
\begin{align}
     Ax - b  = 0, \; A^{T}y + s - c  = 0, \;  XSe  = 0,  \;  (x, \, s)  \ge 0, \label{KKTLP}
\end{align}
where
\begin{align}
    X = diag(x), \;  S = diag(s),  \;  \text{and}\; e = (1, \, \ldots, \, 1)^{T}.
    \label{DIAGXS}
\end{align}
For convenience, we rewrite the optimality condition \eqref{KKTLP}
as the following nonlinear system of equations with nonnegative constraints:
\begin{align}
   F(z) = \begin{bmatrix} Ax - b  \\
                               A^{T}y + s - c \\
                               XS e
  \end{bmatrix} = 0, \; (x, \, s) \ge 0,
  \; \text{and} \; z = (x, \, y, \, s).  \label{NLEQNCON}
\end{align}

\vskip 2mm

It is not difficult to know that the Jacobian matrix $J(z)$ of $F(z)$
has the following form:
\begin{align}
  J(z) = \begin{bmatrix}
    A & 0 & 0 \\
    0 & A^{T} & I \\
    S & 0 & X
  \end{bmatrix}.   \label{JZMATR}
\end{align}
From the third block $XSe = 0$ of equation \eqref{NLEQNCON}, we know that
$x_{i} = 0$ or $s_{i} = 0 \, (i = 1:n) $.  Thus, the Jacobian matrix $J(z)$
of equation \eqref{JZMATR} may be singular, which leads to numerical
difficulties near the solution of the nonlinear system \eqref{NLEQNCON}
for the Newton's method or its variants. In order to overcome this difficulty,
we consider its perturbed system \cite{AG2003,Tanabe1988} as follows:
\begin{align}
  F_{\mu}(z) = F(z) - \begin{bmatrix} 0 \\
                               0 \\
                               \mu e
  \end{bmatrix} = 0, \; (x, \, s) > 0, \; \mu > 0 \;
  \text{and} \; z = (x, \, y, \, s). \label{PNLEX}
\end{align}
The solution $z(\mu)$ of the perturbed system \eqref{PNLEX} defines the
primal-dual central path, and $z(\mu)$ approximates the solution $z^{\ast}$ of
the nonlinear system \eqref{NLEQNCON} when $\mu$ tends to zero
\cite{FMW2007,NW1999,Wright1997,Ye1997}.

\vskip 2mm

We define the strictly feasible region $\mathbb{F}^{0}$ of the problem \eqref{LPND} as
\begin{align}
   \mathbb{F}^{0}  = \left\{(x, \, y, \, s)| Ax = b, \; A^{T}y + s = c,
\; (x, \, s) > 0 \right\}.  \label{SPDFR}
\end{align}
Then, when there is a strictly feasible interior point
$(\bar{x}, \, \bar{y}, \,  \bar{s}) \in \mathbb{F}^{0}$ and the rank of matrix
$A$ is full, the perturbed system \eqref{PNLEX} has a unique solution
(Theorem 2.8, p. 39, \cite{Wright1997}). The existence of its solution
can be derived by the implicit theorem \cite{Doedel2007} and the uniqueness
of its solution can be proved via considering  the strict convexity of
the following penalty problem and the KKT conditions of its optimal solution
 \cite{FM1990}:
\begin{align}
   \min c^{T}x - \mu \sum_{i=1}^{n} \log (x_i) \; \text{subject to}
   \; Ax = b, \label{LOGPENFUN}
\end{align}
where $\mu$ is a positive parameter.

\vskip 2mm

According to the duality theorem of the linear programming (Theorem 13.1, pp. 368-369,
\cite{NW1999}), for any primal-dual feasible solution $(x, \, y, \, s)$, we have
\begin{align}
   c^{T}x \ge  b^{T}y^{\ast} = c^{T}x^{\ast} \ge b^{T}y,  \label{DUATH}
\end{align}
where the triple $(x^{\ast}, \, y^{\ast}, \, s^{\ast})$ is a primal-dual optimal
solution. Moreover, when the positive number $\mu$ is small, the solution
$z^{\ast}(\mu)$ of perturbed system \eqref{PNLEX} is an approximation solution
of nonlinear system \eqref{NLEQNCON}. Consequently, from the duality theorem
\eqref{DUATH}, we know that $x^{\ast}(\mu)$ is an approximation of the optimal
solution of the original linear programming problem \eqref{LPND}. It can be proved
as follows. Since $z^{\ast}(\mu)$ is the primal-dual feasible, from inequality
\eqref{DUATH}, we have
\begin{align}
   c^{T}x^{\ast}(\mu) \ge  b^{T}y^{\ast} = c^{T}x^{\ast} \ge b^{T}y^{\ast}(\mu)
   \label{WDUATH}
\end{align}
and
\begin{align}
   0 \le (x^{\ast}(\mu))^{T}s^{\ast}(\mu)
    = c^{T}x^{\ast}(\mu) -  b^{T}y^{\ast}(\mu). \label{DUALGAP}
\end{align}
From equations \eqref{WDUATH}-\eqref{DUALGAP}, we obtain
\begin{align}
  |c^{T}x^{\ast}(\mu) - c^{T}x^{\ast}| \le c^{T}x^{\ast}(\mu) -  b^{T}y^{\ast}(\mu)
  = (x^{\ast}(\mu))^{T} s^{\ast}(\mu)
  = n \mu.
\end{align}
Therefore, $x^{\ast}(\mu)$ is an approximation of the optimal solution of the
original linear programming problem \eqref{LPND}.

\vskip 2mm

If the damped Newton method is applied to the perturbed system \eqref{PNLEX}
\cite{DS2009,NW1999}, we have
\begin{align}
     z_{k+1} = z_{k} - \alpha_{k} J(z_{k})^{-1}  F_{\mu}(z_{k}), \label{NEWTON}
\end{align}
where $J(z_{k})$ is the Jacobian matrix of $F_{\mu}(z)$. We regard
$z_{k+1} = z(t_{k} + \alpha_{k}), \; z_{k} = z(t_{k})$ and let
$\alpha_{k} \to 0$, then we obtain the continuous Newton flow with the
constraints \cite{AS2015,Branin1972,Davidenko1953,Tanabe1979,LXL2020} of the
perturbed system \eqref{PNLEX}  as follows :
\begin{align}
  \frac{dz(t)}{dt} = - J(z)^{-1}F_{\mu}(z), \hskip 2mm  z = (x, \, y, \, s)
  \; \text{and} \; (x, \, s) > 0.     \label{NEWTONFLOW}
\end{align}
Actually, if we apply an iteration with the explicit Euler method
\cite{SGT2003,YFL1987} for the continuous Newton flow \eqref{NEWTONFLOW}, we
obtain the damped Newton method \eqref{NEWTON}.

\vskip 2mm

Since the Jacobian matrix $J(z) = F'_{\mu}(z)$ may be singular, we reformulate
the continuous Newton flow \eqref{NEWTONFLOW} as the following general
formula \cite{Branin1972,Tanabe1979}:
\begin{align}
    -J(z)\frac{dz(t)}{dt} = F_{\mu}(z), \hskip 2mm  z = (x, \, y, \, s)
    \; \text{and} \; (x, \, s) > 0.  \label{DAEFLOW}
\end{align}
The continuous Newton flow \eqref{DAEFLOW} has some nice properties. We state
one of them as the following property \ref{PRODAEFLOW} \cite{Branin1972,LXL2020,Tanabe1979}.

\vskip 2mm

\begin{property} (Branin \cite{Branin1972} and Tanabe \cite{Tanabe1979})
\label{PRODAEFLOW} Assume that  $z(t)$ is the solution of the continuous Newton flow
\eqref{DAEFLOW}, then $f(z(t)) = \|F_{\mu}(z)\|^{2}$ converges to zero when
$t \to \infty$. That is to say, for every limit point $z^{\ast}$ of $z(t)$, it is
also a solution of the perturbed system \eqref{PNLEX}. Furthermore, every element
$F_{\mu}^{i}(z)$ of $F_{\mu}(z)$ has the same convergence rate $exp(-t)$ and
$z(t)$ can not converge to the solution $z^{\ast}$ of the perturbed system
\eqref{PNLEX} on the finite interval when the initial point $z_{0}$ is not a
solution of the perturbed system \eqref{PNLEX}.
\end{property}
\proof Assume that $z(t)$ is the solution of the continuous Newton flow
\eqref{DAEFLOW}, then we have
\begin{align}
    \frac{d}{dt} \left(exp(t)F_{\mu}(z)\right) = exp(t) J(z) \frac{dz(t)}{dt}
    + exp(t) F_{\mu}(z) = 0. \nonumber
\end{align}
Consequently, we obtain
\begin{align}
     F_{\mu}(z(t)) = F_{\mu}(z_0)exp(-t). \label{FUNPAR}
\end{align}
From equation \eqref{FUNPAR}, it is not difficult to know that every element
$F_{\mu}^{i}(z)$ of $F_{\mu}(z)$ converges to zero with the linear convergence
rate $exp(-t)$ when $t \to \infty$. Thus, if the solution $z(t)$ of the continuous
Newton flow \eqref{DAEFLOW} belongs to a compact set, it has a limit point
$z^{\ast}$ when $t \to \infty$, and this limit point $z^{\ast}$ is
also a solution of the perturbed system \eqref{PNLEX}.

\vskip 2mm

If we assume that the solution $z(t)$ of the continuous Newton flow
\eqref{DAEFLOW} converges to the solution $z^{\ast}$ of the perturbed system
\eqref{PNLEX} on the finite interval $(0, \, T]$, from equation \eqref{FUNPAR},
we have
\begin{align}
     F_{\mu}(z^{\ast}) = F_{\mu}(z_{0}) exp(-T). \label{FLIMT}
\end{align}
Since $z^{\ast}$ is a solution of the perturbed system \eqref{PNLEX}, we have
$F_{\mu}(z^{\ast}) = 0$. By substituting it into equation \eqref{FLIMT}, we
obtain
\begin{align}
     F_{\mu}(z_{0}) = 0. \nonumber
\end{align}
Thus, it contradicts the assumption that $z_{0}$ is not a solution of the perturbed
system \eqref{PNLEX}. Consequently, the solution $z(t)$ of the continuous Newton flow
\eqref{DAEFLOW} can not converge to the solution $z^{\ast}$ of the perturbed system
\eqref{PNLEX} on the finite interval. \qed

\vskip 2mm

\begin{remark} \blue{The inverse $J(x)^{-1}$ of the Jacobian matrix  $J(x)$ can be
regarded as the preconditioner of $F_{\mu}(x)$  such that the every element
$x^{i}(t)$ of $x(t)$ has roughly the same convergence rate and it mitigates the
stiffness of the ODE \eqref{DAEFLOW} \cite{LXL2020}. This property is very
useful since it makes us adopt the explicit ODE method to follow the trajectory
of the Newton flow \eqref{DAEFLOW} efficiently}.
\end{remark}

\vskip 2mm

\subsection{The primal-dual path-following method} \label{SPPFM}

\vskip 2mm

From property \ref{PRODAEFLOW}, we know that the continuous Newton flow
\eqref{DAEFLOW} has the nice global convergence property. However,
when the Jacobian matrix $J(x)$ is singular or nearly singular, the ODE
\eqref{DAEFLOW} is the system of differential-algebraic equations
\cite{AP1998,BCP1996,HW1996} and its trajectory can not be efficiently followed
by the general ODE method such as the backward differentiation formulas (the
built-in subroutine ode15s.m of the MATLAB environment \cite{MATLAB,SGT2003}).
Thus, we need to construct the special method to solve this problem. Furthermore,
we expect that the new method has the global convergence as the homotopy
continuation methods \cite{AG2003,OR2000} and the fast convergence rate as the
traditional optimization methods. In order to achieve these two aims, we consider
the continuation Newton method and the trust-region updating strategy for problem
\eqref{DAEFLOW}.

\vskip 2mm

We apply the implicit Euler method  to the continuous Newton flow \eqref{DAEFLOW}
\cite{AP1998,BCP1996}, then we obtain
\begin{align}
  J(z_{k+1}) \frac{z_{k+1}-z_{k}}{\Delta t_k} = -F_{\mu_{k+1}}(z_{k+1}).
  \label{IMEDAE}
\end{align}
Since the system \eqref{IMEDAE} is nonlinear which is not directly solved, we
seek for its explicit approximation formula. To avoid solving the nonlinear
system of equations, we replace $J(z_{k+1})$ with $J(z_{k})$ and substitute
$F_{\mu_{k+1}}(z_{k+1})$ with its linear approximation
$F_{\mu_{k+1}}(z_{k})+J(z_{k})(z_{k+1}-z_{k})$ into equation \eqref{IMEDAE}.
Then, we obtain a variant of the damped Newton method:
\begin{align}
   z_{k+1}  = z_{k} - \frac{\Delta t_k}{1+\Delta t_k}
   J(z_k)^{-1}F_{\mu_{k+1}}(z_k). \label{SMEDAE}
\end{align}

\vskip 2mm

\begin{remark} If we let $\alpha_{k} = \Delta t_k/(1+\Delta t_k)$ in equation
\eqref{NEWTON}, we obtain the method \eqref{SMEDAE}. However, from the view of
the ODE method, they are different. The damped Newton method \eqref{NEWTON} is
derived from the explicit Euler method applied to the continuous Newton flow
\eqref{DAEFLOW}. Its time-stepping size $\alpha_k$ is restricted by the
numerical stability \cite{HW1996,SGT2003,YFL1987}. That is to say, for the linear
test equation $dx/dt = - \lambda x$, its time-stepping size $\alpha_{k}$ is
restricted by the stable region $|1-\lambda \alpha_{k}| \le 1$. Therefore, the
large time-stepping size $\alpha_{k}$ can not be adopted in the steady-state
phase. The method \eqref{SMEDAE} is derived from \blue{the implicit Euler method
applied to the continuous Newton flow \eqref{DAEFLOW} and the linear approximation
of $F_{\mu_{k+1}}(z_{k+1})$}, and its time-stepping size $\Delta t_k$ is not
restricted by the numerical stability for the linear test equation. Therefore,
the large time-stepping size $\Delta t_{k}$ can be adopted in the steady-state
phase, and the method \eqref{SMEDAE} mimics the Newton method. Consequently,
it has the fast convergence rate near the solution $z^{\ast}$ of the nonlinear
system \eqref{NLEQNCON}. The most of all, the new time-stepping size
$\alpha_{k} = \Delta t_{k}/(\Delta t_{k} + 1)$ is favourable to adopt the
trust-region updating strategy for adaptively adjusting the time-stepping size
$\Delta t_{k}$ such that the continuation method \eqref{SMEDAE} accurately
tracks the trajectory of the continuation Newton flow \eqref{DAEFLOW}
in the transient-state phase and achieves the fast convergence rate
in the steady-state phase.
\end{remark}

\vskip 2mm

We set the parameter $\mu_{k}$ as the average of the residual sum:
\begin{align}
    \mu_{k} =  \frac{\|Ax_{k} - b\|_{1} + \|A^{T}y_{k} + s_{k} - c\|_{1}
    + x_{k}^{T}s_{k}}{n}. \label{UK1DEF}
\end{align}
This selection of $\mu_{k}$ is slightly different to the traditional selection
$\mu_{k} = x_{k}^{T}s_{k}/n$ \cite{FMW2007,NW1999,Wright1997}. According to our
numerical experiments, this selection of $\mu_{k}$ can improve the robustness of
the path-following method. In equation \eqref{SMEDAE}, $\mu_{k+1}$ is
approximated by $\sigma_{k} \mu_{k}$, where the penalty coefficient $\sigma_{k}$ is
simply selected as follows:
\begin{align}
   \sigma_{k} =
   \begin{cases} 0.05, \; \text{when} \;  \mu_{k} > 0.05,  \\
                 \mu_{k}, \; \text{when} \; \mu_{k} \le 0.05.
   \end{cases}  \label{SIGMA}
\end{align}
Thus, from equations \eqref{SMEDAE}-\eqref{SIGMA}, we obtain the following
iteration scheme:
\begin{align}
  \begin{bmatrix}
  A & 0 & 0\\
  0 & A^{T} & I \\
  S_{k} & 0 & X_{k}
 \end{bmatrix}
 \begin{bmatrix} \Delta x_{k} \\
 \Delta y_{k} \\
 \Delta s_{k}
 \end{bmatrix}
 = - \begin{bmatrix}
   Ax_{k} -b \\ A^{T}y_{k} + s_{k} - c \\ X_{k}S_{k}e - \sigma_{k}\mu_{k}e
 \end{bmatrix}
 = - F_{\sigma_{k} \mu_{k}}(z_k) \label{DELTXYSK}
\end{align}
and
\begin{align}
   \left(x_{k+1}, \, y_{k+1}, \, s_{k+1} \right) =
   \left(x_{k}, \, y_{k}, \, s_{k}\right) + \frac{\Delta t_{k}}{1 + \Delta t_{k}}
   \left(\Delta x_{k}, \, \Delta y_{k}, \, \Delta s_{k}\right),
   \label{XYSK1}
\end{align}
where $F_{\mu}(z)$ is defined by equation \eqref{PNLEX}.

\vskip 2mm

When matrix $A$ has full row rank, the linear system \eqref{DELTXYSK} can
be solved by the following three subsystems:
\begin{align}
  AX_{k}S_{k}^{-1}A^{T} \Delta y_{k} & = - \left(r_{p}^{k}
  + AS_{k}^{-1}\left(X_{k}r_{d}^{k} - r_{c}^{k} \right)\right),
  \label{DELTAYK} \\
  \Delta s_{k} & = - r_{d}^{k} - A^{T} \Delta y_{k}, \label{DELTASK} \\
  \Delta x_{k} & = -S_{k}^{-1}\left(X_{k}S_{k}e + X_{k}\Delta s_{k}
  - \sigma_{k} \mu_{k}e\right),   \label{DELTAXK}
\end{align}
where the primal residual $r_{p}^{k}$, the dual residual $r_{d}^{k}$ and the
complementary residual $r_{c}^{k}$ are respectively defined by
\begin{align}
   r_{p}^{k} & = Ax_{k} - b, \label{PRESK} \\
   r_{d}^{k} & = A^{T}y_{k} + s_{k} - c, \label{DRESK} \\
    r_{c}^{d} & =  S_{k}X_{k}e - \sigma_{k} \mu_{k}e. \label{DRESXK}
\end{align}

\vskip 2mm

The matrix $AX_{k}S_{k}^{-1}A^{T}$ becomes very ill-conditioned when
$(x_{k}, \, y_{k}, \, s_{k})$ is close to the solution
$(x^{\ast}, \, y^{\ast}, \, s^{\ast})$ of the nonlinear system \eqref{NLEQNCON}.
Thus, the Cholesky factorization method may fail to solve the linear system
\eqref{DELTAYK} for the large-scale problem. Therefore, we use the QR decomposition
(pp. 247-248, \cite{GV2013}) to solve it as follows:
\begin{align}
     & D_{k}A^{T} = Q_{k}R_{k}, \; D_{k} = \text{diag}(\text{sqrt}(x_{k}./s_{k})),
      \nonumber  \\
     & R_{k}^{T} \Delta y_{k}^{m} = - \left(r_{p}^{k}
     + AS_{k}^{-1}\left(X_{k}r_{d}^{k} - r_{c}^{k} \right)\right), \nonumber \\
     & R_{k} \Delta y_{k} = \Delta y_{k}^{m}.
  \label{DELTAYKQR}
\end{align}
where $Q_{k} \in \Re^{n \times m}$ satisfies $Q_{k}^{T}Q_{k} = I$ and
$R_{k} \in \Re^{m \times m}$ is an upper triangle matrix with full rank.

\vskip 2mm

\subsection{The trust-region updating strategy}

\vskip 2mm

Another issue is how to adaptively adjust the time-stepping size $\Delta t_k$
at every iteration. There is a popular way to control the time-stepping size
based on the trust-region updating strategy
\cite{CGT2000,Deuflhard2004,Higham1999,Luo2010,Luo2012,LXL2020,Yuan2015}. Its
main idea is that the time-stepping size $\Delta t_{k+1}$ will be enlarged
when the linear model $F_{\sigma_{k} \mu_{k}}(z_k) + J(z_k)\Delta z_{k}$ approximates
$F_{\sigma_{k} \mu_{k}}(z_{k}+\Delta z_{k})$ well, and $\Delta t_{k+1}$ will be reduced when
$F_{\sigma_{k} \mu_{k}}(z_k) + J(z_k)\Delta z_{k}$ approximates
$F_{\sigma_{k} \mu_{k}}(z_{k}+J(z_{k})\Delta z_{k})$ badly. We enlarge or reduce the time-stepping size
$\Delta t_k$ at every iteration according to the following ratio:

\begin{align}
  \rho_k = \frac{\|F_{\sigma_{k} \mu_{k}}(z_{k})\|-\|F_{\sigma_{k} \mu_{k}}(z_{k}+\Delta z_{k})\|}
    {\|F_{\sigma_{k} \mu_{k}}(z_{k})\|
    - \|F_{\sigma_{k} \mu_{k}}(z_{k})+J(z_{k})\Delta z_{k}\|}. \label{RHOK}
\end{align}
A particular adjustment strategy is given as follows:
\begin{align}
   \Delta t_{k+1} =
     \begin{cases}
          2 \Delta t_k, \; &{\text{if} \; 0 \leq \left|1- \rho_k \right| \le \eta_1
           \; \text{and} \; (x_{k+1}, \, s_{k+1}) > 0,} \\
         \Delta t_k, \; &{\text{else if} \; \eta_1 < \left|1 - \rho_k \right| < \eta_2
         \; \text{and} \; (x_{k+1}, \, s_{k+1}) > 0, }\\
    \frac{1}{2}  \Delta t_k, \; &{\text{others}, }
    \end{cases} \label{TSK1}
\end{align}
where the constants are selected as $\eta_1 = 0.25, \; \eta_2 = 0.75$,
according to our numerical experiments. When $\rho_{k} \ge \eta_{a}$ and
$(x_{k+1}, \, s_{k+1}) > 0$, we accept the trial step
and set
\begin{align}
      (x_{k+1}, \, y_{k+1}, \, s_{k+1})  = (x_{k}, \, y_{k}, \, s_{k})
      + \frac{\Delta t_{k}}{1+\Delta t_{k}} (\Delta x_{k}, \, \Delta y_{k}, \, \Delta s_{k}),
      \label{ACCPXK1}
\end{align}
where $\eta_{a}$ is a small positive number such as $\eta_{a} = 1.0\times 10^{-6}$.
Otherwise, we discard the trial step and set
\begin{align}
     (x_{k+1}, \, y_{k+1}, \, s_{k+1})  = (x_{k}, \, y_{k}, \, s_{k}). \label{NOACXK1}
\end{align}

\vskip 2mm

\begin{remark}
This new time-stepping size selection based on the trust-region updating strategy
has some advantages compared to the traditional line search strategy.
If we use the line search strategy and the damped Newton method \eqref{NEWTON}
to track the trajectory $z(t)$ of the continuous Newton flow \eqref{DAEFLOW},
in order to achieve the fast convergence rate in the steady-state phase,
the time-stepping size $\alpha_{k}$ of the damped Newton method is tried from 1
and reduced by the half with many times at every iteration. Since the linear
model $F_{\sigma_{k} \mu_{k}}(z_{k}) + J(z_{k})\Delta z_{k}$ may not
approximate $F_{\sigma_{k}\mu_{k}}(z_{k}+\Delta z_{k})$ well in the transient-state phase,
the time-stepping size $\alpha_{k}$ will be small. Consequently, the line search
strategy consumes the unnecessary trial steps in the transient-state phase.
However, the selection of the time-stepping size $\Delta t_{k}$ based on the
trust-region updating strategy \eqref{RHOK}-\eqref{TSK1} can overcome this shortcoming.
\end{remark}

\subsection{The treatment of rank-deficient problems}

\vskip 2mm

For a real-world problem, the rank of matrix $A$ may be deficient and
the constraints are even inconsistent when the right-hand-side vector $b$
has small noise \cite{LLS2020}. However, the constraints of the original problem
are intrinsically consistent. For the rank-deficient problem with the consistent
constraints, there are some efficiently pre-solving methods to eliminate the
redundant constraints in references \cite{AA1995,Andersen1995,MS2003}. Here,
in order to handle the inconsistent system of constraints, we consider the
following least-squares approximation problem:
\begin{align}
   \min_{x \in \Re^{n}} \; \|Ax - b\|^{2}. \label{LLSP}
\end{align}
Then, by solving problem \eqref{LLSP}, we obtain the consistent system of
constraints.

\vskip 2mm

Firstly, we use the QR factorization with column pivoting (pp. 276-278,
\cite{GV2013}) to factor $A$ into a product of an orthogonal matrix
$Q \in \Re^{m \times m}$ and an upper triangular matrix $R \in \Re^{m \times n}$:
\begin{align}
      AP = QR = \begin{bmatrix} Q_{1} | Q_{2} \end{bmatrix}
      \begin{bmatrix} R_{1} \\ 0     \end{bmatrix}
      = Q_{1}R_{1}, \label{ATQR}
\end{align}
where $r = \text{rank}(A), \; Q_{1} = Q(1:m, \, 1:r), \; R_{1} = R(1:r, \, 1:n)$,
and $P \in \Re^{n \times n}$ is a permutation matrix. Then, from equation \eqref{ATQR},
we know that problem \eqref{LLSP} equals the following problem
\begin{align}
      \min_{x \in \Re^{n}} \; \left\| R_{1}P^{T}x - Q_{1}^{T}b\right\|^{2}. \label{BAPROB}
\end{align}
By solving problem \eqref{BAPROB}, we obtain its solution $x$ as follows:
\begin{align}
    R_{1}\tilde{x} = b_{r}, \; x = P\tilde{x}, \label{APPCON}
\end{align}
where $b_{r} = Q_{1}^{T}b$.

\vskip 2mm

Therefore, when the constraints of problem \eqref{LPND} are consistent, problem
\eqref{LPND} equals the following linear programming problem:
\begin{align}
     \min_{x \in \Re^{n}} \; c_{r}^{T} \tilde{x}
      \; \text{subject to} \; R_{1} \tilde{x} = b_{r}, \; \tilde{x} \ge 0,
      \label{LOLPND}
\end{align}
where $b_r = Q_{1}^{T}b$ and $c_{r} = P^{T}c$.
When the constraints of problem \eqref{LPND} are inconsistent, the constraints
of problem \eqref{LOLPND} are the least-squares approximation of constraints of
problem \eqref{LPND}. Consequently, in subsection \ref{SPPFM}, we replace matrix
$A$, vector $b$ and vector $c$  with matrix $R_{1}$, vector $b_{r}$ and
vector $c_{r}$, respectively, then the primal-dual path-following method
can handle the rank-deficient problem.

\vskip 2mm

From the reduced linear programming problem \eqref{LOLPND}, we obtain its KKT
conditions as follows:
\begin{align}
     R_{1}\tilde{x} - b_{r} & = 0, \label{LOPFCON} \\
     R_{1}^{T} \tilde{y} + \tilde{s} - c_{r} & = 0, \label{LODFCON} \\
     \tilde{X}\tilde{S} e & = 0, \; i = 1, \, 2, \, \ldots, n, \label{LOCCON} \\
     (\tilde{x}, \, \tilde{s}) & \ge 0, \label{LONNCON}
\end{align}
where $\tilde{X} = diag(\tilde{x}), \; \tilde{S} = diag(\tilde{s})$, $b_{r} = Q_{1}^{T}b$
and $c_{r} = P^{T}c$. Thus, from the QR decomposition \eqref{ATQR} of matrix $A$
and equation \eqref{KKTLP}, we can recover the solution $(x, \, y, \, s)$ of
equation \eqref{KKTLP} as follows:
\begin{align}
      x = P\tilde{x}, \; y = Q_{1}\tilde{y}, \; s = P\tilde{s}, \label{RSOL}
\end{align}
where $(\tilde{x}, \, \tilde{y}, \, \tilde{s})$ is a solution of
equations \eqref{LOPFCON}-\eqref{LONNCON}.

\vskip 2mm

\begin{remark}
\blue{The preprocessing strategies of the redundant elimination in references
\cite{AA1995,Andersen1995,MS2003} are for empty rows and columns,
the row or column singletons, duplicate rows or columns, forcing and dominated
constraints, and finding the linear dependency based on the Gaussian
elimination. Some of those techniques such as for empty rows and columns,
the row or column singletons, duplicate rows or columns can be also applied to the
inconsistent system. However, those preprocessing strategies in
references \cite{AA1995,Andersen1995,MS2003} can not transform an inconsistent
system to a consistent system. Thus, in order to handle the inconsistent system,
we can replace the QR decomposition with the Gaussian elimination method  after
the preprocessing strategies such as for empty rows or columns, the row or column
singletons, duplicate rows or colums. Here, for simplicity, we directly use the QR
decomposition as the preprocessing strategy and do not use the preprocessing
strategies in references \cite{AA1995,Andersen1995,MS2003}.}
\end{remark}

According to the above discussions, we give the detailed descriptions of the
primal-dual path-following method and the trust-region updating strategy for
linear the programming problem \eqref{LPND} in Algorithm \ref{PNMTRLP}.

\vskip 2mm

\begin{algorithm}
    \renewcommand{\algorithmicrequire}{\textbf{Input:}}
	\renewcommand{\algorithmicensure}{\textbf{Output:}}
	\caption{Primal-dual path-following methods and the trust-region updating strategy for
    linear programming (The PFMTRLP method)}
    \label{PNMTRLP}
    \begin{algorithmic}[1]
        \REQUIRE ~~\\
          matrix $A \in \Re^{m \times n}$, vectors $b \in \Re^{m}$ and $c \in \Re^{n}$
        for the problem:
        $\min_{x \in \Re^{n}} c^{T}x \; \text{subject to} \; Ax  = b, \; x \ge 0$.
		\ENSURE ~~ \\
        the primal-dual optimal solution: $(solx, \, soly, \, sols)$.   \\
        the maximum error of KKT conditions: $KKTError = \max \{\|A*solx -b\|_{\infty}, \,
        \|A^{T}*soly + sols -c \|_{\infty}, \, \|solx.*soly\|_{\infty}\}$.
        \STATE Initialize parameters: $\eta_{a} = 10^{-6}, \; \eta_1 = 0.25,
        \; \eta_2 = 0.75, \; \epsilon = 10^{-6}, \; \Delta t_0 = 0.9,
        \; \text{maxit} = 100$.
        \STATE Factor matrix $A$ by using the QR decomposition \eqref{ATQR}.
        Set $b_{r} = Q_{1}^{T}b$ and $c_{r} = P^{T}c$.
        \STATE Set $\text{bigM} = \max(\max(\text{abs}(R_{1})))$; $\text{bigM}
        = \max(\|b_{r}\|_{\infty}, \, \|c_{r}\|_{\infty}, \,  \text{bigM})$.

        \STATE Initialize $x_0$ = bigMfac*bigM*ones(n, 1), $s_0 = x_0$, $y_0$ = zeros(r, 1).
        \STATE flag\_success\_trialstep = 1, $\text{itc} = 0, \; k = 0$.
        \WHILE {(itc $<$ maxit)}
           \IF{(flag\_success\_trialstep == 1)}
              \STATE Set itc = itc + 1.
              \STATE  Compute $F_{k}^{1} = R_{1}x_k - b_{r}$, $F_{k}^{2} = R_{1}^{T}y_k + s_k - c_{r}$,
              $F_{k}^{3} = X_{k}s_k$, $\mu_{k} = (\|F_{k}^{1}\|_{1} + \|F_{k}^{2}\|_{1} + \|F_{k}^{3}\|_{1})/n$.
              \STATE  Compute $\text{Resk} = \|[F_{k}^{1}; F_{k}^{2}; F_{k}^{3}]\|_{\infty}$.
              \IF{($\text{Resk} < \epsilon$)}
                 \STATE break;
              \ENDIF
              \STATE Set $\sigma_k = \min(0.05, \, \mu_k)$.
              \STATE Compute $F_{k}^{3} = F_{k}^{3} - \sigma_{k}\mu_{k}*\text{ones}(n,1)$.
              Set $F_k = [F_{k}^{1}; F_{k}^{2}; F_{k}^{3}]$. \\
              \STATE By using the QR decomposition to solve \eqref{DELTAYKQR}, we obtain $\Delta y_{k}$.
              Compute $\Delta s_{k}   = - F_{k}^{2}- Q_{1}\Delta y_{k}$, and
                $\Delta x_{k}   = -S_{k}^{-1}\left(F_{k}^{3} + X_{k}\Delta s_{k} \right).
              $
           \ENDIF
           \STATE Set $(x_{k+1}, \, y_{k+1}, \, s_{k+1}) = (x_{k}, \, y_{k}, \, s_{k})
              + \frac{\Delta t_{k}}{1 + \Delta t_{k}} (\Delta x_{k}, \, \Delta y_{k}, \, \Delta s_{k})$.
           \STATE Compute $F_{k+1}^{1} = R_{1}^{T}x_{k+1} - b_{r}$, $F_{k+1}^{2} = R_{1}^{T}y_{k+1} + s_{k+1} - c_{r}$,
               $F_{k+1}^{3} = X_{k+1}s_{k+1} - \sigma_{k}\mu_{k}*\text{ones}(n,1)$,
              $F_{k+1} = [F_{k+1}^{1};F_{k+1}^{2};F_{k+1}^{3}]$.
           \STATE Compute $LinAppF_{k+1} = [F_{k+1}^{1}; F_{k+1}^{2};(F_{k}^{3} + X_{k}(s_{k+1}-s_{k}) + S_{k}(x_{k+1}-x_{k}))]$.
           \STATE Compute the ratio $\rho_{k} = (\|F_{k}\|-\|F_{k+1}\|)/(\|F_{k}\|-\|LinAppF_{k+1}\|)$.
           \IF{($(|\rho_{k} - 1| \le \eta_{1}) \&\& (x_{k+1}, s_{k+1}) > 0)$)}
              \STATE  Set $\Delta t_{k+1} = 2\Delta t_{k}$;
           \ELSIF{($(\eta_{1} < |\rho_{k} - 1| \le \eta_{2}) \&\& (x_{k+1}, s_{k+1}) > 0))$}
               \STATE Set $\Delta t_{k+1} = \Delta t_{k}$;
           \ELSE
               \STATE Set $\Delta t_{k+1} = 0.5\Delta t_{k}$;
           \ENDIF
           \IF{$((\rho_{k} \ge \eta_{a}) \&\& (x_{k+1}, s_{k+1}) > 0))$}
               \STATE Accept the trial point $(x_{k+1}, \, y_{k+1}, \, s_{k+1})$. Set
               flag\_success\_trialstep = 1.
           \ELSE
               \STATE Set $(x_{k+1}, \, y_{k+1}, \, s_{k+1}) = (x_{k}, \, y_{k}, s_{k})$,
               flag\_success\_trialstep = 0.
           \ENDIF
           \STATE Set $k \leftarrow k+1$.
        \ENDWHILE
        \STATE Return $(solx, \, soly, \, sols) =  \left(Px_{k}, \,
        Q_{1}y_{k}, \, Ps_{k}\right)$, $KKTError = \|F_{k}\|_{\infty}$.
   \end{algorithmic}
\end{algorithm}

\vskip 2mm

\section{Algorithm analysis}

\vskip 2mm

We define the one-sided neighborhood $\mathbb{N}_{-\infty}(\gamma)$ as
\begin{align}
     \mathbb{N}_{-\infty}(\gamma) =
     \{(x, \, y, \, s)  \in \mathbb{F}^{0} | XSe \ge \gamma \mu e \}, \label{OSINFN}
\end{align}
where $X = \text{diag}(x), \; S = \text{diag}(s), \; e = \text{ones}(n,\,1), \;
\mu = x^{T}s/n$ and $\gamma$ is a small positive constant such as $\gamma = 10^{-3}$.
In order to simplify the convergence analysis of Algorithm \ref{PNMTRLP}, we assume
that (i) the initial point $(x_{0}, \, s_{0})$ is strictly primal-dual feasible,
and (ii) the time-stepping size $\Delta t_{k}$ is selected such that
$(x_{k+1}, \, y_{k+1}, \, s_{k+1}) \in \mathbb{N}_{-\infty}(\gamma)$. Without the loss
of generality, we assume that the row rank of matrix $A\in \Re^{m\times n}$ is full.

\vskip 2mm

\begin{lemma} \label{LEMPROC}
Assume $(x_{k}, \, y_{k}, \, s_{k}) \in \mathbb{N}_{-\infty}(\gamma)$, then
there exists a sufficiently small positive number $\delta_{k}$ such that
$(x_{k}(\alpha), \, y_{k}(\alpha), \, s_{k}(\alpha)) \in \mathbb{N}_{-\infty}(\gamma)$
holds when $0 < \alpha \le \delta_{k}$, where
$(x_{k}(\alpha), \, y_{k}(\alpha), \, s_{k}(\alpha))$ is defined by
\begin{align}
     (x_{k}(\alpha), \, y_{k}(\alpha), \, s_{k}(\alpha))
     = (x_{k}, \, y_{k}, \, s_{k})
     + \alpha (\Delta x_{k}, \, \Delta y_{k}, \, \Delta s_{k}),
      \label{LINST}
\end{align}
and $(\Delta x_{k}, \, \Delta y_{k}, \, \Delta s_{k})$ is the solution of the linear
system \eqref{DELTXYSK}.
\end{lemma}

\vskip 2mm

\proof Since $(x_{k}, \, y_{k}, \, s_{k})$ is a primal-dual feasible point,
from equation \eqref{DELTXYSK}, we obtain
\begin{align}
    A \Delta x_{k} = 0, \;  A^{T} \Delta y_{k} + \Delta s_{k} = 0.
    \label{NULLSP}
\end{align}
Consequently, from equations \eqref{LINST}-\eqref{NULLSP}, we have
\begin{align}
    \Delta x_{k}^{T} \Delta s_{k} = - \Delta x_{k}^{T}A^{T} \Delta y_{k}
    = (A \Delta x_{k})^{T} \Delta y_{k} = 0,    \label{DELSXSUMZ} \\
     Ax_{k}(\alpha) = b, \; A^{T}y_{k}(\alpha) + s_{k}(\alpha) = c.
     \label{PDFR}
\end{align}
From equation \eqref{LINST}, we have
\begin{align}
     & X_{k}(\alpha)S_{k}(\alpha)  =
     (X_{k} + \alpha \Delta X_{k})(S_{k} + \alpha \Delta S_{k}) \nonumber \\
     & \quad = X_{k}S_{k} + \alpha  (X_{k}\Delta S_{k} + S_{k}\Delta X_{k})
      + \alpha^{2} \Delta X_{k} \Delta S_{k}. \label{XKSKPRO}
\end{align}
By replacing $X_{k}\Delta S_{k} + S_{k}\Delta X_{k}$ with equation \eqref{DELTAXK}
into equation \eqref{XKSKPRO}, we obtain
\begin{align}
    & X_{k}(\alpha)S_{k}(\alpha)e = X_{k}S_{k}e +
     \alpha (\sigma_{k}\mu_{k}e - X_{k}S_{k}e) + \alpha^{2} \Delta X_{k} \Delta S_{k}e
    \nonumber \\
    & \quad = (1 - \alpha)X_{k}S_{k}e + \alpha \sigma_{k}\mu_{k}e
    + \alpha^{2} \Delta X_{k} \Delta S_{k}e. \label{XKSKVEC}
\end{align}
From equations \eqref{DELSXSUMZ}-\eqref{XKSKVEC}, we have
\begin{align}
  \mu_{k}(\alpha) =  \frac{1}{n} e^{T}X_{k}(\alpha)S_{k}(\alpha)e
  = (1 - (1-\sigma_{k})\alpha) \mu_{k}. \label{MUALPHA}
\end{align}

We denote
\begin{align}
   \beta_{max}^{k}
     = \max_{1 \le i \le n} \{|\Delta x_{k}^{i}|, \, |\Delta s_{k}^{i}|\}.
    \label{BETACON}
\end{align}
Then, from equation \eqref{XKSKVEC}, we have
\begin{align}
    X_{k}(\alpha)S_{k}(\alpha)e \ge \left((1 - \alpha)\gamma \mu_{k}
    + \alpha \sigma_{k}\mu_{k} - \alpha^{2} (\beta_{max}^{k})^{2}\right)e.
    \label{XKSKGEL}
\end{align}
From equations \eqref{MUALPHA} and \eqref{XKSKGEL}, we know that the
proximity condition
\begin{align}
    X_{k}(\alpha)S_{k}(\alpha) e \ge \gamma \mu_{k}(\alpha) \nonumber
\end{align}
holds, provided that
\begin{align}
    (1 - \alpha)\gamma \mu_{k} + \alpha \sigma_{k}\mu_{k}
     - \alpha^{2}(\beta_{max}^{k})^{2}
     \ge \gamma (1- \alpha + \alpha \sigma_{k})\mu_{k}. \nonumber
\end{align}
By reformulating the above expression, we obtain
\begin{align}
      \alpha (1-\gamma) \sigma_{k} \mu_{k} \ge \alpha^{2} (\beta_{max}^{k})^{2}.
    \label{ALPKGE}
\end{align}
We choose
\begin{align}
     \delta_{k} = \frac{(1-\gamma)\sigma_{k}\mu_{k}}{(\beta_{max}^{k})^{2}}.
     \label{ALPKLE}
\end{align}
Then, inequality \eqref{ALPKGE} is true when $0< \alpha \le \delta_{k}$.  \qed

\vskip 2mm

In the following Lemma \ref{LEMXSBOUN}, we give the lower bounded estimation
of $(x_{k}, \, s_{k})$.

\vskip 2mm

\begin{lemma} \label{LEMXSBOUN}
Assume that $(x_{0}, \, y_{0}, \, s_{0}) > 0$ is a primal-dual feasible
point and $(x_{k}, \, y_{k}, \, s_{k})$ $(k = 0, \, 1, \, \ldots)$ generated
by Algorithm \ref{PNMTRLP} satisfy the proximity condition \eqref{OSINFN}.
Furthermore, if there exists a constant $C_{\mu}$ such that
\begin{align}
   \mu_{k} \ge C_{\mu} > 0 \label{MUKGEC}
\end{align}
holds for all $k = 0, \, 1, \, 2, \, \ldots$, there exist two positive constants
$C_{min}$ and $C_{max}$ such that
\begin{align}
    0 < C_{min} \le \min_{1\le i \le n} \{x_{k}^{i}, \; s_{k}^{i}\}
    \le \max_{1 \le i \le n}  \{x_{k}^{i},\; s_{k}^{i} \}
    \le C_{max}    \label{BOUNDXS}
\end{align}
holds for all $k = 0, \, 1, \, 2, \, \ldots$.
\end{lemma}

\vskip 2mm

\proof Since $(x_{k}, \, y_{k}, s_{k})$ is generated by Algorithm \ref{PNMTRLP}
and $(x_{k}, \, y_{k}, s_{k})$ is a primal-dual feasible point, from equation
\eqref{MUALPHA}, we have
\begin{align}
   \mu_{k+1} = \frac{x_{k+1}^{T}s_{k+1}}{n}
    = (1 - (1-\sigma_{k})\alpha_{k}) \mu_{k} \le \mu_{k},
    \; k = 0, \, 1, \, \ldots. \nonumber
\end{align}
Consequently, we obtain
\begin{align}
    \mu_{k+1} = \prod_{i = 0}^{k} (1 - (1-\sigma_{i})\alpha_{i}) \mu_{i}
    \le \mu_{0}, \; k = 0, \, 1, \, \ldots. \label{UKUPBOUN}
\end{align}

\vskip 2mm

From equation \eqref{PDFR}, we have
\begin{align}
     A (x_{k} - x_{0}) = 0, \; A^{T}(y_{k} - y_{0}) + (s_{k} - s_{0}) = 0.
\end{align}
Consequently, we obtain
\begin{align}
    (x_{k} - x_{0})^{T}(s_{k} - s_{0}) = - (x_{k} - x_{0})^{T} A^{T}(y_{k} - y_{0})
    = 0. \nonumber
\end{align}
By rearranging this expression and using the property \eqref{UKUPBOUN}, we obtain
\begin{align}
    x_{k}^{T}s_{0} + s_{k}^{T}x_{0}
     = x_{k}^{T}s_{k} + x_{0}^{T}s_{0} \le n (\mu_{k} + \mu_{0})
    \le 2n \mu_{0}. \nonumber
\end{align}
Consequently, we obtain
\begin{align}
   x_{k}^{i} \le  \frac{2n \mu_{0}} {\min_{1 \le j \le n} \{s_{0}^{j}\}}\;
   \text{and} \;
   s_{k}^{i} \le \frac{2n \mu_{0}}{\min_{1 \le j \le n} \{x_{0}^{j}\}}, \;
   1 \le i \le n, \; k = 0,\, 1, \, \ldots. \nonumber
\end{align}
Therefore, if we select
\begin{align}
    C_{max} = \max \left\{\frac{2n \mu_{0}}{\min_{1 \le j \le n} \{s_{0}^{j}\}},
     \; \frac{2n \mu_{0}}{\min_{1 \le j \le n} \{x_{0}^{j}\}}\right\}, \nonumber
\end{align}
 we obtain
\begin{align}
    \max_{1 \le i \le n} \{x_{k}^{i}, s_{k}^{i}\} \le C_{max}, \;
    k = 0, \, 1, \, \ldots.
    \label{XSKUPBOUN}
\end{align}

\vskip 2mm

On the other hand, from the assumption \eqref{MUKGEC} and the proximity
condition \eqref{OSINFN}, we have
\begin{align}
  x_{k}^{i}s_{k}^{i}  \ge \gamma \mu_{k} \ge \gamma C_{\mu}, \; 1 \le i \le n,
  \; k = 0, \, 1, \, \ldots.  \nonumber
\end{align}
By combining it with the estimation \eqref{XSKUPBOUN} of $(x_{k}, \, s_{k})$, we
obtain
\begin{align}
    x_{k}^{i} \ge \frac{\gamma C_{\mu}}{\max_{1 \le j \le n}\{s_{k}^{j}\}}
    \ge  \frac{\gamma C_{\mu}}{C_{max}},  \; \text{and}\;
     s_{k}^{i} \ge \frac{\gamma C_{\mu}}{\max_{1 \le j \le n} \{x_{k}^{j}\}}
    \ge  \frac{\gamma C_{\mu}}{C_{max}}, \; k = 0, \, 1, \, \ldots. \label{XSLBOUD}
\end{align}
We select $C_{min} = \gamma C_{\mu}/C_{max}$.  Then, from equation \eqref{XSLBOUD},
we obtain
\begin{align}
    \min_{1 \le i \le n} \{x_{k}^{i}, \, s_{k}^{i}\} \ge C_{min}, \;
    k = 0, \, 1, \, 2, \, \ldots.  \nonumber
\end{align}
 \qed

\vskip 2mm

\begin{lemma} \label{LEMDELSXUB}
Assume that $(x_{0}, \, y_{0}, \, s_{0}) > 0$ is a primal-dual feasible
point and $(x_{k}, \, y_{k}, \, s_{k})$ \, $(k = 0, \, 1, \, \ldots)$ generated
by Algorithm \ref{PNMTRLP} satisfy the proximity condition \eqref{OSINFN}.
Furthermore, if the assumption \eqref{MUKGEC} holds for all
$k = 0, \, 1, \, \ldots$, there exit two positive constants $C_{\Delta x}$
and $C_{\Delta s}$ such that
\begin{align}
   \|\Delta s_{k}\| \le C_{\Delta s} \; \text{and} \; \|\Delta x_{k}\| \le C_{\Delta x}
   \label{DELXSUB}
\end{align}
hold for all $k = 0, \, 1, \, \ldots$.
\end{lemma}

\proof Factorize matrix $A$ with the singular value decomposition (pp. 76-80, \cite{GV2013}):
\begin{align}
     A = U\Sigma V^{T}, \; \Sigma =
     \begin{bmatrix} \Sigma_{r} & 0 \\ 0 &  0 \end{bmatrix}, \;
     \Sigma_{r} = \text{diag}(\lambda_{1}, \, \lambda_{2}, \, \ldots, \, \lambda_{r})
     \succ 0,      \label{SVDA}
\end{align}
where $U \in \Re^{m \times m}$ and $V \in \Re^{n \times n}$ are orthogonal
matrices, and the rank of matrix $A$ equals $r$. Then, from the bounded estimation
\eqref{BOUNDXS} of $(x_{k}, \, s_{k})$, we have
\begin{align}
     z AX_{k}S^{-1}_{k}A^{T}z \ge \frac{C_{min}}{C_{max}} \|Az\|^2
     \ge \frac{C_{min}\lambda_{min}^{2}}{C_{max}} \|z\|^{2}, \;
     k = 0, \, 1, \, \ldots, \; \forall z \in \Re^{n},     \label{AXSINLB}
\end{align}
and
\begin{align}
    z AX_{k}S_{k}^{-1}A^{T}z \le \frac{C_{max}}{C_{min}} \|Az\|^2
    \le \frac{C_{max}\lambda_{max}^{2}}{C_{min}} \|z\|^{2}, \;
    k = 0, \, 1, \ldots, \;  \forall z \in \Re^{n},     \label{AXSINUB}
\end{align}
where $\lambda_{min}$ and $\lambda_{max}$ are the smallest and largest
singular values of matrix $A$, respectively.

\vskip 2mm

From equations \eqref{DELTAYK}, \eqref{BOUNDXS} and
\eqref{AXSINLB}-\eqref{AXSINUB}, we obtain
\begin{align}
    & \frac{C_{min}\lambda_{min}^{2}}{C_{max}} \|\Delta y_{k} \|^{2}
     \le \Delta y_{k}^{T} \left(AX_{k}S_{k}^{-1}A^{T}\right) \Delta y_{k}
      = \Delta y_{k}^{T} AS_{k}^{-1}(X_{k}S_{k}e - \sigma_{k}\mu_{k}e)
     \nonumber \\
  & \quad \le \|\Delta y_{k}\| \|A\| \left\|S_{k}^{-1}\right\|
  \left\|X_{k}S_{k}e - \sigma_{k}\mu_{k}e\right\|
   \le  \|\Delta y_{k}\|  \frac{\lambda_{max}}{C_{min}}
  (\|X_{k}S_{k}e\| + n \sigma_{k}\mu_{k}) \nonumber \\
   & \quad \le \|\Delta y_{k}\|  \frac{\lambda_{max}}{C_{min}}
  (\|X_{k}S_{k}e\|_{1} + n \sigma_{k}\mu_{k})
  = \|\Delta y_{k}\|  \frac{\lambda_{max}}{C_{min}} (1+\sigma_{k})n \mu_{k}.
  \nonumber
\end{align}
That is to say, we obtain
\begin{align}
    \|\Delta y_{k}\|
     \le \frac{C_{max}\lambda_{max}}{C_{min}^{2}\lambda_{min}^{2}}(1+\sigma_{k}) n\mu_{k}
    \le \frac{C_{max}\lambda_{max}}{C_{min}^{2}\lambda_{min}^{2}}2n\mu_{k}
     \le \frac{C_{max}\lambda_{max}}{C_{min}^{2}\lambda_{min}^{2}}2n\mu_{0}.
    \label{DELYKUB}
\end{align}
The second inequality of equation \eqref{DELYKUB} can be inferred by
$\sigma_{k} \le 1$ from equation \eqref{SIGMA}. The last inequality of equation
\eqref{DELYKUB} can be inferred by  $\mu_{k} \le \mu_{0}$ from equation
\eqref{UKUPBOUN}. Therefore, from equation \eqref{DELTASK} and equation
\eqref{DELYKUB}, we have
\begin{align}
    \|\Delta s_{k}\| = \|-A^{T}\Delta y_{k}\| \le \|A^{T}\| \|\Delta y_{k} \|
    \le \frac{C_{max}\lambda^{2}_{max}}{C_{min}^{2}\lambda_{min}^{2}}2n\mu_{0}.
    \label{DELSKUPBD}
\end{align}
We set $C_{\Delta s} = (C_{max}\lambda^{2}_{max}2n\mu_{0})/(C_{min}^{2}\lambda_{min}^{2})$.
Thus, from equation \eqref{DELSKUPBD}, we prove the first part of equation
\eqref{DELXSUB}.

\vskip 2mm

From equations \eqref{DELTAXK},  \eqref{BOUNDXS} and the first part of
equation \eqref{DELXSUB}, we have
\begin{align}
    & \|\Delta x_{k}\|   = \left\|-S_{k}^{-1}(X_{k}S_{k}e + X_{k}\Delta s_{k}
     - \sigma \mu_{k}e)\right\| \le \left\|S_{k}^{-1}\right\|
    \|X_{k}S_{k}e + X_{k}\Delta s_{k} - \sigma_{k} \mu_{k}e\| \nonumber \\
    & \quad \le \frac{1}{C_{min}}
     \left(\|X_{k}S_{k}e\|+\|X_{k}\Delta s_{k}\|+\|\sigma_{k}\mu_{k}e \|\right)
     \nonumber \\
    & \quad \le \frac{1}{C_{min}}\left(\|X_{k}S_{k}e\|_{1}
     + \|X_{k}\| \|\Delta s_{k}\|+ n \sigma_{k}\mu_{k} \|\right)
     \nonumber \\
     & \quad \le \frac{1}{C_{min}}
    \left(n \mu_{k} + C_{max} C_{\Delta s} + n \sigma_{k} \mu_{k}\right)
      \le \frac{1}{C_{min}} \left(2n\mu_{0} + C_{max} C_{\Delta s}\right).
        \label{DELXKUPBD}
\end{align}
The last inequality of equation \eqref{DELXKUPBD} can be inferred by
$\sigma_{k} \le 1$ from equation \eqref{SIGMA} and $\mu_{k} \le \mu_{0}$
from equation \eqref{UKUPBOUN}. We set
$C_{\Delta x} = \left(2n\mu_{0} + C_{max} C_{\Delta s}\right)/C_{min}$.
Thus, from equation \eqref{DELXKUPBD}, we also prove the second part of
equation \eqref{DELXSUB}.  \qed

\vskip 2mm

\begin{lemma} \label{LEMDELTLB}
Assume that $(x_{0}, \, y_{0}, \, s_{0}) > 0$ is a primal-dual feasible
point and $(x_{k}, \, y_{k}, \, s_{k})$ \, $(k = 0, \, 1, \, \ldots)$ generated
by Algorithm \ref{PNMTRLP} satisfy the proximity condition \eqref{OSINFN}.
Furthermore, if the assumption \eqref{MUKGEC} holds for all
$k = 0, \, 1, \, \ldots$, there exits a positive constant $C_{\Delta t}$ such that
\begin{align}
    \Delta t_{k} \ge C_{\Delta t} > 0 \label{DELTLB}
\end{align}
holds for all $k = 0, \, 1, \, 2, \ldots$.
\end{lemma}

\vskip 2mm

\proof From equations \eqref{RHOK}, \eqref{DELTXYSK}-\eqref{XYSK1}, we have
\begin{align}
   |\rho_k - 1| & = \left|\frac{\|F_{\sigma_{k} \mu_{k}}(z_{k})\|
   -\|F_{\sigma_{k} \mu_{k}}(z_{k}+\Delta z_{k})\|}{\|F_{\sigma_{k} \mu_{k}}(z_{k})\|
    - \|F_{\sigma_{k} \mu_{k}}(z_{k})+J(z_{k}) \Delta z_{k}\|} - 1\right| \nonumber \\
    & = \left|\frac{\|F_{\sigma_{k} \mu_{k}}(z_{k})\|-\|F_{\sigma_{k} \mu_{k}}(z_{k}+\Delta z_{k})\|}
    {\|F_{\sigma_{k} \mu_{k}}(z_{k})\| - \|F_{\sigma_{k} \mu_{k}}(z_{k}) -
    (\Delta t_{k})/(1+\Delta t_{k})F_{\sigma_{k}\mu_{k}}(z_{k})\|} - 1\right| \nonumber \\
    & = \frac{ \left|\|F_{\sigma_{k}\mu_{k}}(z_{k})\|
    - (1+\Delta t_{k})\|F_{\sigma_{k}\mu_{k}}(z_{k}+\Delta z_{k})\|\right|}
    {\Delta t_{k} \|F_{\sigma_{k}\mu_{k}}(z_{k})\|} \nonumber \\
    & \le \frac{ \left\|(1+\Delta t_{k})\left(F_{\sigma_{k}\mu_{k}}(z_{k}+\Delta z_{k})
      - F_{\sigma_{k}\mu_{k}}(z_{k})\right)
      + \Delta t_{k} F_{\sigma_{k}\mu_{k}}(z_{k})\right\|}
    {\Delta t_{k} \|F_{\sigma_{k}\mu_{k}}(z_{k})\|} \nonumber \\
    & = \frac{\Delta t_{k}}{1+\Delta t_{k}} \frac{\|\Delta X_{k}\Delta S_{k}e\|}
    {\|X_{k}s_{k} - \sigma_{k}\mu_{k}e\|}.
     \label{RHOMINUS1}
\end{align}
The last equality of equation \eqref{RHOMINUS1} can be inferred from
\begin{align}
    & F_{\sigma_{k}\mu_{k}}(z_{k}+\Delta z_{k}) - F_{\sigma_{k}\mu_{k}}(z_{k}) 
    = J(z_{k})\Delta z_{k} + 
    \left(\frac{\Delta t_{k}}{1+\Delta t_{k}}\right)^{2}\Delta X_{k}\Delta S_{k}e \nonumber \\
   & \quad = - \frac{\Delta t_{k}}{1+\Delta t_{k}} F_{\sigma_{k}\mu_{k}}(z_{k})
   + \left(\frac{\Delta t_{k}}{1+\Delta t_{k}}\right)^{2}\Delta X_{k}\Delta S_{k}e \nonumber.
\end{align}

\vskip 2mm

On the other hand, from the property $\|a\| \ge a_{i} \, (i = 1, \, 2, \ldots, n)$,
we have
\begin{align}
    \|X_{k}s_{k} - \sigma_{k}\mu_{k}e\|  \ge x_{k}^{i}s_{k}^{i} - \sigma_{k}\mu_{k},
    \; i = 1, \, 2, \, \ldots, n.   \nonumber
\end{align}
By summing the $n$ components of the above two sides and $\mu_{k} = x_{k}^{T}s_{k}/n$,
we obtain
\begin{align}
    \|X_{k}s_{k} - \sigma_{k}\mu_{k}e\| \ge (1-\sigma_{k}) \mu_{k}
    \ge (1 - 0.05) C_{\mu} = 0.95C_{\mu}. \label{XSLBD}
\end{align}
The second inequality of equation \eqref{XSLBD} can be inferred by
$\sigma_{k} \le 0.05$ from equation \eqref{SIGMA} and the assumption
$\mu_{k} \ge C_{\mu}$ from equation \eqref{MUKGEC}.

\vskip 2mm

Thus, from the bounded estimation \eqref{DELXSUB} of
$(\Delta x_{k}, \, \Delta s_{k})$ and equations \eqref{RHOMINUS1}-\eqref{XSLBD},
we obtain
\begin{align}
  |\rho_{k} - 1| & \le \frac{\Delta t_{k}}{1+\Delta t_{k}}
  \frac{\|\Delta X_{k} \Delta S_{k}e\|} {0.95C_{\mu}}
  \le \frac{\Delta t_{k}}{1+\Delta t_{k}}
    \frac{\|\Delta X_{k}\| \|\Delta s_{k}\|} {0.95C_{\mu}} \nonumber \\
  & \le \frac{\Delta t_{k}}{1+\Delta t_{k}}
    \frac{\|\Delta x_{k}\| \|\Delta s_{k}\|} {0.95C_{\mu}}
   \le \frac{\Delta t_{k}}{1+\Delta t_{k}}
    \frac{C_{\Delta x} C_{\Delta s}} {0.95C_{\mu}}
   \le \eta_{1}, \nonumber
\end{align}
provided that
\begin{align}
    \Delta t_{k} \le \frac{0.95 C_{\mu}\eta_{1}}
    {C_{\Delta x}C_{\Delta s}}. \label{DLETKLECON}
\end{align}
Thus, if we assume that $K$ is the first index such that
$\Delta t_{K}$ satisfies equation \eqref{DLETKLECON}, according to the
trust-region updating formula \eqref{TSK1}, $\Delta t_{K+1}$ will be
enlarged. Therefore, we prove the result \eqref{DELTLB} if we set
\begin{align}
    C_{\Delta t}  = \frac{0.95 C_{\mu}\eta_{1}}
    {2C_{\Delta x}C_{\Delta s}}. \nonumber
\end{align}
  \qed

\vskip 2mm

According to the above discussions, we give the global convergence analysis
of Algorithm \ref{PNMTRLP}.

\begin{theorem}  Assume that $(x_{0}, \, y_{0}, \, s_{0}) > 0$ is a primal-dual
feasible point and $(x_{k}, \, y_{k}, \, s_{k})$ \, $(k = 0, \, 1, \, \ldots)$
generated by Algorithm \ref{PNMTRLP} satisfy the proximity condition \eqref{OSINFN}.
Then, we have
\begin{align}
    \lim_{k \to \infty} \mu_{k} = 0, \label{MUKTOZ}
\end{align}
where $\mu_{k} = x_{k}^{T}s_{k}/n$.
\end{theorem}

\proof Assume that there exists  a positive constant $C_{\mu}$ such that
\begin{align}
   \mu_{k} \ge C_{\mu} \label{UKGELB}
\end{align}
holds for all $k = 0, \, 1, \, \ldots$. Then, according to the result of Lemma
\ref{LEMDELTLB}, we know that there exists a positive constant $C_{\Delta t}$
such that $\Delta t_{k} \ge C_{\Delta t}$ holds for all $k = 0, \, 1, \ldots$.
Therefore, from equation \eqref{MUALPHA}, we have
\begin{align}
   \mu_{k+1} & = \mu_{k}(\alpha_{k}) = (1 - (1 - \sigma_{k})\alpha_{k}) \mu_{k}
   \le \left(1 - (1-\sigma_{k}) \frac{C_{\Delta t}}{1+C_{\Delta t}}\right) \mu_{k}
   \nonumber \\
   & \le \left(1 - 0.95 \frac{C_{\Delta t}}{1+C_{\Delta t}}\right) \mu_{k}
   \le \left(1 - 0.95 \frac{C_{\Delta t}}{1+C_{\Delta t}}\right)^{k+1} \mu_{0}.
   \label{UK1LE}
\end{align}
The second inequality of equation \eqref{UK1LE} can be inferred by
$\sigma_{k} \le 0.05$ from equation \eqref{SIGMA}. Thus, we have $\mu_{k} \to 0$,
 which contradicts the assumption \eqref{UKGELB}.
Consequently, we obtain $\lim_{k \to \infty} \inf{\mu_{k}} = 0$. Since $\mu_{k}$
is monotonically decreasing, it is not difficult to know
$\lim_{ k \to \infty} \mu_{k} = 0$. Furthermore, we obtain
$\lim_{k \to \infty} \|X_{k}s_{k}\| = 0$ from
$\|X_{k}s_{k}\| \le \|X_{k}s_{k}\|_{1} = n \mu_{k}$ and $(x_{k}, s_{k}) > 0$.  \qed

\vskip 2mm

\begin{remark} Our analysis framework of Algorithm \ref{PNMTRLP} is the same as
that of the classical primal-dual interior method (pp. 411-413, \cite{NW1999}).
However, the result of Lemma \ref{LEMDELTLB} is new.
\end{remark}

\section{Numerical experiments}

\vskip 2mm

In this section, we test Algorithm \ref{PNMTRLP} (the PFMTRLP method) for some
linear programming problems with full rank matrices or rank-deficient matrices,
and compare it with the traditional path-following method
(pathfollow.m, p. 210, \cite{FMW2007}) and the state-of-the-art predictor-corrector
algorithm (the built-in subroutine linprog.m of the MATLAB environment
\cite{MATLAB,Mehrotra1992,Zhang1998}).

\vskip 2mm

The tolerable errors of three methods are all set by $\epsilon = 10^{-6}$. We
use the maximum absolute error (KKTError) of the KKT condition
\eqref{KKTLP} and the primal-dual gap $x^{T}s$ to measure the
error between the numerical optimal solution and the theoretical optimal
solution.

\vskip 2mm

\subsection{The problem with full rank}

For the standard linear programming problem with full rank,
the sparse matrix $A$ of given density 0.2 is randomly generated and we
choose feasible $x, \, y, \, s$ at random, with $x$ and $s$ each about
half-full. The dimension of matrix $A$ varies from $10\times100$ to
$300\times3000$. One of its implementation is given by Algorithm
\ref{FULLMATPRO} (p. 210, \cite{FMW2007}). According to Algorithm
\ref{FULLMATPRO}, we randomly generate 30 standard linear programming problems
with full rank matrices.

\vskip 2mm

\begin{algorithm}
    \renewcommand{\algorithmicrequire}{\textbf{Input:}}
	\renewcommand{\algorithmicensure}{\textbf{Output:}}
	\caption{The standard linear programming problem with full rank}
    \label{FULLMATPRO}
    \begin{algorithmic}[1]
    \REQUIRE ~~\\
       the number of equality constraints: $m$; \\
       the number of unknown variables: $n$.
    \ENSURE ~~ \\
       matrix $A$ and vectors $b \in \Re^{m}$ and $c \in \Re^{n}$.
    \STATE density=0.2;
    \STATE A = sprandn(m,n,density); \% Generate a sparse matrix of give density.
    \STATE xfeas = [rand(n/2,1); zeros(n-(n/2),1)];
    \STATE sfeas = [zeros(n/2,1); rand(n-(n/2),1)];
    \STATE xfeas = xfeas(randperm(n));
    \STATE sfeas = sfeas(randperm(n));
    \STATE yfeas = (rand(m,1)-0.5)*4;\\
    \% Choose b and c to make this (x,y,s) feasible.
    \STATE b = A*xfeas;
    \STATE c = A'*yfeas + sfeas;
    \end{algorithmic}
\end{algorithm}

\vskip 2mm

\begin{figure}[htbp]
  \centering
  \begin{minipage}[t]{0.49\linewidth}
     \centering
     \subfigure[The number of iterations]{
        \includegraphics[width=1\textwidth, height=0.25\textheight]{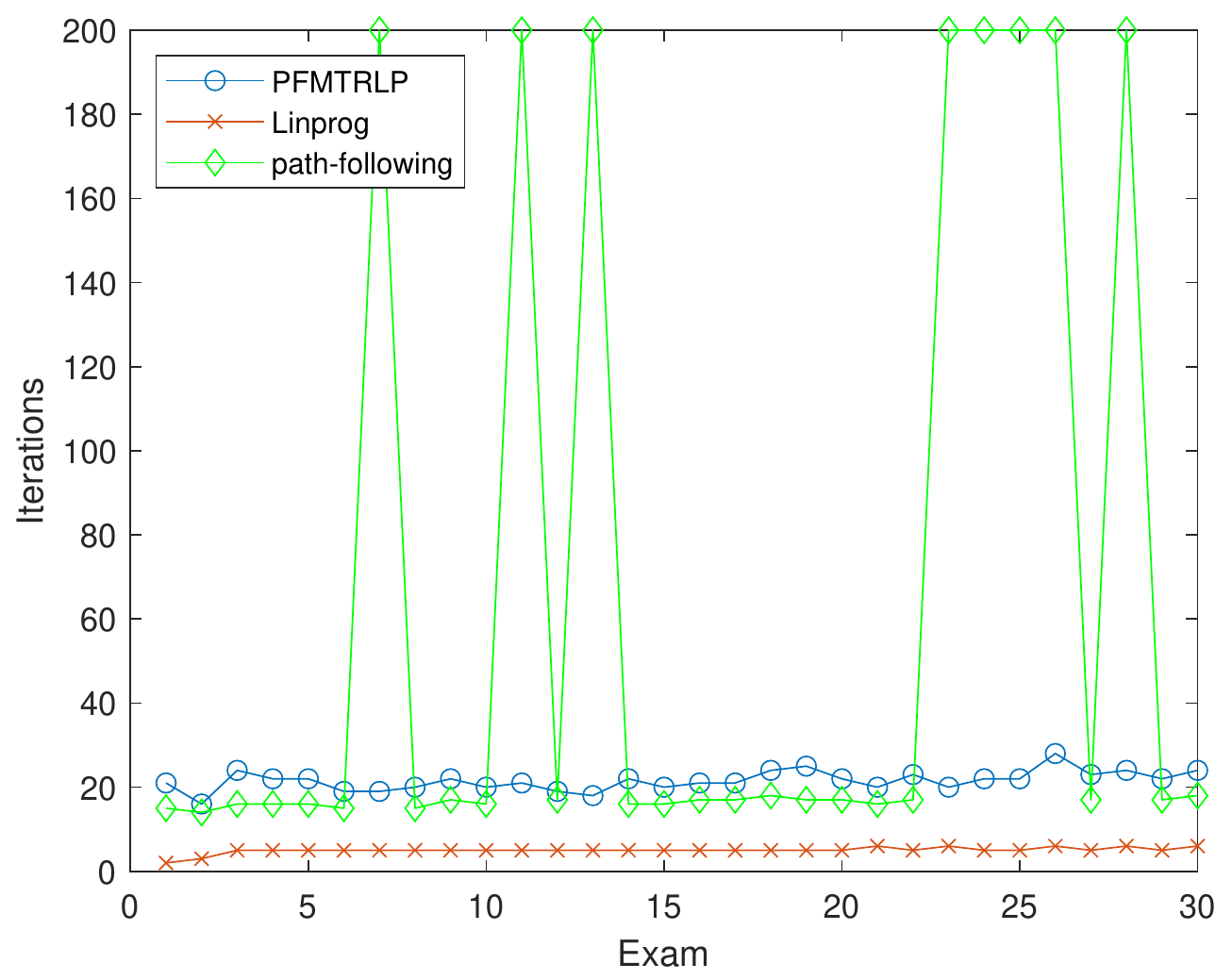}
     }
  \end{minipage}
  \begin{minipage}[t]{0.49\linewidth}
     \centering
     \subfigure[The computational time]{
         \includegraphics[width=1\textwidth,height=0.25\textheight]{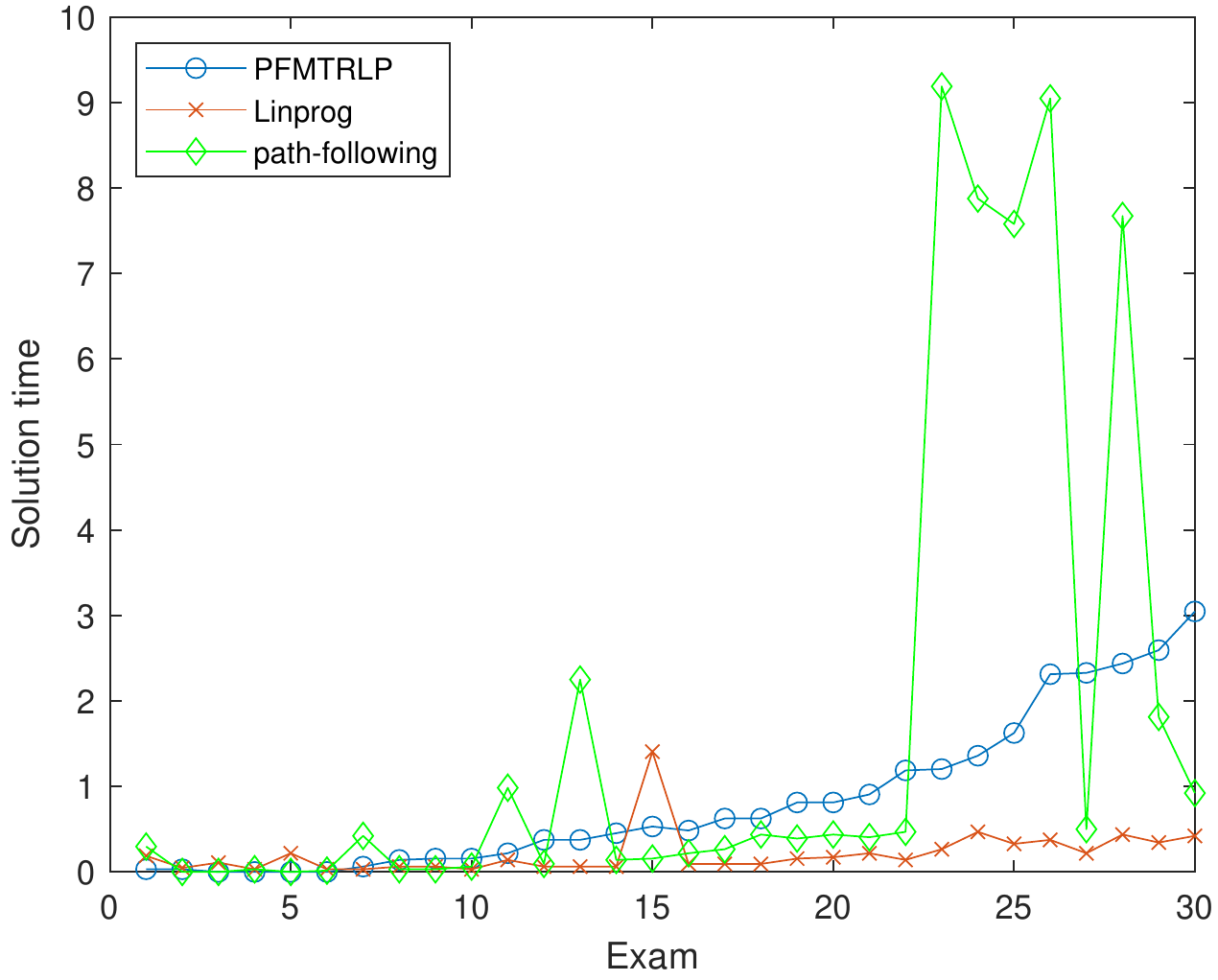}
     }
  \end{minipage}
  \caption{The number of iterations and the computational time.} %
  \label{FIGFULLRANK}
\end{figure}

For those 30 test problems, we compare Algorithm \ref{PNMTRLP}
(the PFMTRLP method), Mehrotra's predictor-corrector algorithm (the
subroutine linprog.m of the MATLAB environment), and the path-following method
(the subroutine pathfollow.m). The numerical results are arranged in Table
\ref{TABFULLRANK} and illustrated in Figure \ref{FIGFULLRANK}. The left
sub-figure of Figure \ref{FIGFULLRANK} represents the number of iterations
and the right sub-figure represents the consumed CPU time. From Table
\ref{TABFULLRANK}, we find that PFMTRLP and linprog.m can solve all test problems,
and their \emph{KKTErrors} are small. However, pathfollow.m  cannot
perform well for some higher-dimensional problems, such as exam
$7, \, 11, \, 13,  \, 23, \, 24, \, 25, \, 26, \, 28$,
since their solutions do not satisfy the KKT conditions. From Figure
\ref{FIGFULLRANK}, we also find that linprog.m performs the best,
and the number of its iterations is less than 20. The number of iterations of
PFMTRLP is around 20, and the number of iterations of pathfollow.m often
reaches the maximum number (i.e. $200$ iterations). Therefore, PFMTRLP
is also an efficient and robust path-following method for
the linear programming problem with full rank.

\begin{table}[htbp]
  \centering
  \fontsize{7}{7}\selectfont
  \caption{Numerical results of problems with full rank matrices.}
  \label{TABFULLRANK}
    \begin{tabular}{|l|l|l|l|l|l|l|}
    \hline
    \multirow{2}{*}{Problem ($m \times n$, $r$)} &
    \multicolumn{2}{c|}{PFMTRLP} & \multicolumn{2}{c|}{linprog} & \multicolumn{2}{c|}{pathfollow} \cr \cline{2-7}
           & KKTError  & Gap      & KKTError & Gap      & KKTError & Gap      \cr\hline
    Exam. 1 ($10\times100$, 10)     & 3.77E-06 & 3.61E-05 & 1.98E-07    & 1.18E-08  & 1.03E-07  & 6.27E-06  \cr\hline
    Exam. 2 ($20\times200$, 20)     & 2.07E-06 & 1.46E-04 & 3.20E-10    & 8.83E-14  & 6.62E-08  & 1.04E-05  \cr\hline
    Exam. 3 ($30\times300$, 30)     & 9.68E-06 & 1.28E-03 & 1.07E-11    & 3.10E-10  & 1.83E-09  & 4.48E-07  \cr\hline
    Exam. 4 ($40\times400$, 40)     & 3.16E-06 & 5.67E-04 & 6.55E-08    & 1.15E-06  & 1.34E-07  & 4.42E-05  \cr\hline
    Exam. 5 ($50\times500$, 50)     & 1.14E-05 & 2.47E-03 & 5.12E-07    & 3.79E-05  & 1.25E-08  & 4.55E-06  \cr\hline
    Exam. 6 ($60\times600$, 60)     & 1.42E-06 & 2.26E-04 & 6.83E-09    & 2.90E-07  & 3.95E-09  & 1.79E-06  \cr\hline
    Exam. 7 ($70\times3700$, 70)    & 1.05E-04 & 1.27E-02 & 4.21E-07    & 3.71E-05  & \red{4.61E+04}  & 4.48E-04  \cr\hline
    Exam. 8 ($80\times800$, 80)     & 7.67E-06 & 2.36E-03 & 4.40E-09    & 5.09E-07  & 4.03E-07  & 2.64E-04  \cr\hline
    Exam. 9 ($90\times900$, 90)     & 7.67E-06 & 2.43E-03 & 1.14E-09    & 1.08E-07  & 2.62E-08  & 1.64E-05  \cr\hline
    Exam. 10 ($100\times1000$, 100) & 1.93E-05 & 7.86E-03 & 2.34E-12    & 3.67E-11  & 9.09E-09  & 6.65E-06  \cr\hline
    Exam. 11 ($1100\times1100$, 110)& 9.00E-05 & 3.54E-02 & 6.36E-08    & 1.13E-05  & \red{8.52E+04}  & 9.33E-04  \cr\hline
    Exam. 12 ($120\times1200$, 120) & 1.46E-05 & 5.08E-03 & 3.84E-08    & 1.18E-05  & 1.08E-09  & 4.78E-07  \cr\hline
    Exam. 13 ($130\times1300$, 130) & 8.44E-06 & 3.78E-03 & 1.37E-10    & 2.67E-08  & \red{1.23E+05}  & 2.17E-03  \cr\hline
    Exam. 14 ($140\times1400$, 140) & 3.33E-06 & 1.55E-03 & 3.53E-07    & 3.78E-05  & 8.73E-07  & 8.96E-04  \cr\hline
    Exam. 15 ($150\times1500$, 150) & 7.23E-05 & 2.71E-02 & 3.59E-07    & 4.16E-05  & 1.25E-06  & 1.40E-03  \cr\hline
    Exam. 16 ($160\times1600$, 160) & 1.19E-05 & 6.79E-03 & 3.58E-08    & 1.22E-05  & 3.07E-07  & 4.06E-04  \cr\hline
    Exam. 17 ($170\times1700$, 170) & 4.45E-05 & 2.25E-02 & 6.33E-11    & 1.37E-08  & 3.53E-09  & 4.74E-06  \cr\hline
    Exam. 18 ($180\times1800$, 180) & 1.22E-04 & 6.32E-02 & 8.85E-07    & 3.28E-04  & 1.88E-08  & 2.37E-05  \cr\hline
    Exam. 19 ($190\times1900$, 190) & 6.51E-05 & 5.81E-02 & 8.86E-08    & 6.25E-06  & 2.06E-07  & 3.24E-04  \cr\hline
    Exam. 20 ($200\times2000$, 200) & 2.49E-05 & 1.61E-02 & 7.56E-07    & 4.67E-04  & 1.23E-06  & 1.91E-03  \cr\hline
    Exam. 21 ($210\times2100$, 210) & 1.42E-04 & 7.20E-02 & 3.50E-13    & 7.61E-11  & 1.33E-06  & 1.92E-03  \cr\hline
    Exam. 22 ($220\times2200$, 220) & 1.59E-05 & 1.04E-02 & 1.64E-07    & 7.08E-05  & 2.72E-08  & 4.33E-05  \cr\hline
    Exam. 23 ($230\times2300$, 230) & 3.64E-04 & 2.28E-01 & 8.65E-14    & 2.62E-11  & \red{2.75E+05}  & 1.48E-02  \cr\hline
    Exam. 24 ($240\times24000$, 240)& 9.80E-05 & 8.36E-02 & 6.76E-07    & 2.19E-04  & \red{2.39E+05}  & 2.19E-02  \cr\hline
    Exam. 25 ($250\times2500$, 250) & 3.06E-04 & 1.80E-01 & 8.40E-08    & 3.59E-05  & \red{2.13E+05}  & 8.76E-02  \cr\hline
    Exam. 26 ($260\times2600$, 260) & 1.21E-05 & 8.87E-03 & 7.98E-09    & 2.15E-07  & \red{5.59E+05}  & 3.12E-01  \cr\hline
    Exam. 27 ($270\times2700$, 2700)& 1.15E-04 & 1.44E-01 & 4.79E-07    & 1.33E-04  & 2.06E-08  & 4.20E-05   \cr\hline
    Exam. 28 ($280\times2800$, 280) & 3.33E-05 & 4.27E-02 & 2.58E-13    & 1.13E-10  & \red{4.94E+05}  & 1.87E-02  \cr\hline
    Exam. 29 ($290\times2900$, 290) & 4.42E-05 & 4.38E-02 & 3.31E-07    & 1.04E-04  & 3.81E-08  & 8.19E-05  \cr\hline
    Exam. 30 ($300\times3000$, 300) & 7.84E-05 & 4.21E-02 & 1.82E-08    & 1.24E-05  & 5.38E-08  & 1.30E-04   \cr\hline
    \end{tabular}%
\end{table}%

\vskip 2mm

\vskip 2mm

\subsection{The rank-deficient problem with noisy data}

\vskip 2mm

For a real-world problem, the rank of matrix $A$ in problem \eqref{LPND}
may be deficient and the constraints are even inconsistent when the
right-hand-side vector $b$ has small noise. However, the constraints of the
original problem are intrinsically consistent. In order to evaluate the effect
of PFMTRLP handling those problems, we select some rank-deficient problems from
the NETLIB collection \cite{NETLIB} as test problems and compare it with
linprog.m for those problems with or without the small noisy data.

\vskip 2mm

The numerical results of the problems without the noisy data are arranged
in Table \ref{NETLIBORIGINAL}. Then, we set $b = b + \text{rand}(m,\, 1)*\epsilon$
for those test problems, where $\epsilon = 10^{-5}$. The numerical results of the
problems with the small noise are arranged in Table \ref{NETLIBNOISE}. From
Tables \ref{NETLIBORIGINAL} and \ref{NETLIBNOISE}, we can find that PFMTRLP can
solve all those problems with or without the small noise.  However, from Table
\ref{NETLIBNOISE}, we find that linprog.m can not solve some problems with the
small noise since inprog.m outputs $NaN$ for those problems. Furthermore, although
linprog.m can solve some problems, the KKT errors or the primal-dual gaps are
large. For those problems, we conclude that linprog.m also fails to solve them.
Therefore, from Table \ref{NETLIBNOISE}, we find that PFMTRLP is more robust than
linprog.m for the rank-deficient problem with the small noisy data.

\vskip 2mm

\begin{table}[htbp]
  \newcommand{\tabincell}[2]{\begin{tabular}{@{}#1@{}}#2\end{tabular}}
  \scriptsize
  \centering
  \fontsize{7}{7}\selectfont
  \caption{Numerical results of some rank-deficient problems from NETLIB.}
  \label{NETLIBORIGINAL}
    \begin{tabular}{|l|l|l|l|l|l|l|l|l|}
    \hline
    \multirow{2}{*}{Problem ($m \times n$, $r$)} &
    \multicolumn{4}{c|}{PFMTRLP} & \multicolumn{4}{c|}{linprog}  \cr \cline{2-9}
           & KKTError  & Gap   & Iter  & CPU   & KKTError & Gap      & Iter & CPU      \cr\hline
    \tabincell{c}{lp\_brandy \\ ($220 \times 303$, 193)}
            & 2.37E-03 & 3.57E-02 & 38    & 0.19    & 2.44E-08 & 1.14E-08 & 16    & 0.09    \cr\hline
    \tabincell{c}{lp\_bore3d \\  ($233 \times 334$, 231)}        & 3.00E-08 & 3.32E-06 & 43    & 0.22    & 5.67E-10 & \red{1.37E+03} & 17    & 0.24    \cr\hline
    \tabincell{c}{lp\_wood1p \\ ($244 \times 2595$, 243)}       & 5.09E-05 & 3.09E-05 & 100   & 5.28    & 1.80E-12 & 3.16E-11 & 20    & 0.26    \cr\hline
    \tabincell{c}{lp\_scorpion \\ ($388 \times 466$, 358)}      & 6.16E-07 & 8.65E-05 & 37    & 0.64    & 2.12E-13 & 5.08E-08 & 14    & 0.05    \cr\hline
    \tabincell{c}{lp\_ship04s \\ ($402 \times 1506$, 360)}      & 7.06E-07 & 2.71E-03 & 36    & 1.83    & 2.39E-07 & 2.24E-04 & 13    & 0.03    \cr\hline
    \tabincell{c}{lp\_ship04l \\ ($402 \times 2166$, 360)}      & 1.79E-07 & 1.00E-03 & 36    & 2.48    & 1.04E-11 & 1.02E-04 & 12    & 0.08    \cr\hline
    \tabincell{c}{lp\_degen2 \\ ($444 \times 757$, 442)}        & 8.90E-07 & 4.45E-05 & 35    & 1.23    & 2.79E-13 & 9.43E-10 & 14    & 0.40    \cr\hline
    \tabincell{c}{lp\_bnl1 \\ ($643 \times 1586$, 642)}       & 7.40E-07 & 6.50E-04 & 97    & 11.33   & 3.36E-09 & 3.16E-06 & 26    & 0.08    \cr\hline
    \tabincell{c}{lp\_ship08s \\ ($778 \times 2467$, 712)}      & 4.34E-07 & 3.29E-03 & 41    & 12.72   & 2.12E-11 & 1.98E-07 & 14    & 0.04    \cr\hline
    \tabincell{c}{lp\_qap8 \\ ($912 \times 1632$, 742)}         & 6.84E-07 & 1.82E-04 & 22    & 4.78    & 2.67E-11 & 6.48E-07 & 9     & 0.49    \cr\hline
    \tabincell{c}{lp\_25fv47 \\ ($821 \times 1876$, 820)}       & 5.74E-07 & 9.20E-04 & 80    & 14.67   & 3.35E-10 & 6.10E-11 & 25    & 0.24    \cr\hline
    \tabincell{c}{lp\_ship08l \\ ($778 \times 4363$, 712)}      & 9.95E-07 & 9.87E-03 & 40    & 18.84   & 1.07E-09 & 4.28E-03 & 15    & 0.15    \cr\hline
    \tabincell{c}{lp\_ship12l \\ ($1151 \times 5533$, 1041)}    & 9.22E-07 & 4.21E-03 & 47    & 37.06   & 4.08E-10 & 1.15E-02 & 15    & 0.06    \cr\hline
    \tabincell{c}{lp\_ship12s \\ ($1151 \times 2869$, 1042)}    & 3.01E-07 & 9.36E-04 & 47    & 26.11   & 1.75E-11 & 1.50E-05 & 16    & 0.05    \cr\hline
    \tabincell{c}{lp\_degen3 \\ ($1503 \times 2604$, 1501)}     & 4.79E-07 & 1.01E-05 & 48    & 27.61   & 5.90E-10 & 3.55E-08 & 21    & 0.52    \cr\hline
    \tabincell{c}{lp\_qap12 \\ ($3192 \times 8856$, 2794)}     & 7.00E-07 & 6.86E-04 & 25    & 228.75  & 1.01E-05 & 2.95E-04 & 85    & 218.06  \cr\hline
    \tabincell{c}{lp\_cre\_c \\ ($3068 \times 6411$, 2981)}     & 7.97E-07 & 4.79E-01 & 90    & 516.67  & 5.15E-10 & 4.50E-02 & 27    & 0.21    \cr\hline
    \tabincell{c}{lp\_cre\_a \\ ($3516 \times 7248$, 3423)}    & 7.86E-07 & 7.02E-01 & 83    & 707.98  & 7.26E-09 & 7.10E-02 & 28    & 0.24    \cr\hline
    \end{tabular}%
\end{table}%

\begin{table}[htbp]
  \newcommand{\tabincell}[2]{\begin{tabular}{@{}#1@{}}#2\end{tabular}}
  \scriptsize
  \centering
  \fontsize{7}{7}\selectfont
  \caption{Numerical results of rank-deficient problems with noise $\epsilon = 10^{-5}$.}
  \label{NETLIBNOISE}
    \begin{tabular}{|l|l|l|l|l|l|l|l|l|}
    \hline
    \multirow{2}{*}{Problem ($m \times n$, $r$)} &
    \multicolumn{4}{c|}{PFMTRLP} & \multicolumn{4}{c|}{linprog}  \cr \cline{2-9}
           & KKTError  & Gap   & Iter  & CPU   & KKTError & Gap      & Iter & CPU      \cr\hline
    \tabincell{c}{lp\_brandy \\ ($220 \times 303$, 193)}        & 4.21E-02 & 2.63E+00 & 28    & 0.06    & \red{Failed}    & Failed    & 0     & 0.01  \cr\hline
    \tabincell{c}{lp\_bore3d \\ ($233 \times 334$, 231)}        & 1.84E-01 & 1.69E+01 & 42    & 0.19    & \red{Failed}    & Failed    & 0     & 0.01  \cr\hline
    \tabincell{c}{lp\_scorpion \\ ($388 \times 466$, 358)}      & 5.03E-02 & 1.05E+01 & 18    & 0.36    & \red{Failed}    & Failed    & 0     & 0.01  \cr\hline
    \tabincell{c}{lp\_degen2 \\ ($444 \times 757$, 442)}        & 1.41E-03 & 9.63E-02 & 28    & 1.22    & 2.92E-04  & \red{2.12E+10}  & 52    & 0.37  \cr\hline
    \tabincell{c}{lp\_ship04s \\ ($402 \times 1506$, 360)}     & 3.91E-02 & \red{2.93E+02} & 23    & 1.56    & \red{Failed}    & Failed    & 0     & 0.00  \cr\hline
    \tabincell{c}{lp\_bnl1 \\ ($643 \times 1586$, 642)}        & 3.49E-03 & 2.68E+00 & 61    & 7.70    & \red{Failed}    & Failed    & 0     & 0.01  \cr\hline
    \tabincell{c}{lp\_qap8 \\ ($912 \times 1632$, 742)}         & 8.08E-07 & 2.19E-04 & 22    & 4.58    & 6.24E-04  & \red{6.93E+02}  & 6     & 0.35  \cr\hline
    \tabincell{c}{lp\_25fv47 \\ ($821 \times 1876$, 820)}       & 6.34E-07 & 1.01E-03 & 80    & 14.41   & \red{Failed}    & Failed    & 0     & 0.01  \cr\hline
    \tabincell{c}{lp\_ship04l \\ ($402 \times 2166$, 360)}      & 9.44E-03 & 7.35E+01 & 25    & 2.03    & \red{Failed}    & Failed    & 0     & 0.00  \cr\hline
    \tabincell{c}{lp\_ship08s \\ ($778 \times 2467$, 712)}      & 1.69E-02 & \red{1.72E+02} & 24    & 8.31    & \red{Failed}    & Failed    & 0     & 0.00  \cr\hline
    \tabincell{c}{lp\_wood1p \\ ($244 \times 2595$, 243)}       & 8.91E-04 & 6.99E-04 & 60    & 3.16    & 2.41E-04  & \red{4.64E+07}  & 33    & 0.62  \cr\hline
    \tabincell{c}{lp\_degen3 \\ ($1503 \times 2604$, 1501)}     & 9.56E-04 & 2.92E-02 & 37    & 24.06   & 6.84E-04  & \red{1.89E+05}  & 77    & 2.53  \cr\hline
    \tabincell{c}{lp\_ship12s \\ ($1151 \times 2869$, 1042)}    & 1.36E-02 & 6.26E+01 & 31    & 18.39   & \red{Failed}    & Failed    & 0     & 0.03  \cr\hline
    \tabincell{c}{lp\_ship08l \\ ($778 \times 4363$, 712)}      & 1.36E-02 & \red{1.55E+02} & 24    & 12.80   & \red{Failed}    & Failed    & 0     & 0.00  \cr\hline
    \tabincell{c}{lp\_ship12l \\ ($1151 \times 5533$, 1041)}    & 1.18E-02 & 5.88E+01 & 31    & 29.44   & \red{Failed}    & Failed    & 0     & 0.00  \cr\hline
    \tabincell{c}{lp\_cre\_c \\ ($3068 \times 6411$, 2981)}     & 5.81E-02 & \red{3.84E+04} & 43    & 277.59  & \red{Failed}    & Failed    & 0     & 0.01  \cr\hline
    \tabincell{c}{lp\_cre\_a \\ ($3516 \times 7248$, 3423)}     & 7.37E-02 & \red{6.73E+04} & 37    & 350.22  & \red{Failed}    & Failed    & 0     & 0.00  \cr\hline
    \tabincell{c}{lp\_qap12 \\ ($3192 \times 8856$, 2794)}      & 7.01E-07 & 6.86E-04 & 25    & 240.84  & 2.07E-02  & \red{1.40E+03}  & 10    & 24.80 \cr\hline
    \end{tabular}%
\end{table}%

\section{Conclusions}

\vskip 2mm

For the rank-deficient linear programming problem, we give
a preprocessing method based on the QR decomposition with column pivoting.
Then, we consider the primal-dual path-following and the trust-region updating
strategy for the postprocessing problem. Finally, we prove that the global
convergence of the new method when the initial point is strictly prima-dual
feasible. According to our numerical experiments, the new method (PFMTRLP)
is more robust than the path-following methods such as pathfollow.m (p. 210,
\cite{FMW2007}) and linprog.m \cite{MATLAB,Mehrotra1992,Zhang1998} for
the rank-deficient problem with the small noisy data. Therefore, PFMTRLP
is worth exploring further as a primal-dual path-following method with
the new adaptive step size selection based on the trust-region updating
stratgy. The computational efficiency of PFMTRLP has a room to improve.

\vskip 2mm

\section*{Acknowledgments}
This work was supported in part by Grant 61876199 from National Natural Science
Foundation of China, Grant YBWL2011085 from Huawei Technologies
Co., Ltd., and Grant YJCB2011003HI from the Innovation Research Program of Huawei
Technologies Co., Ltd.. The first author is grateful to professor Li-zhi Liao for
introducing him the interior-point methods when he visited Hongkong Baptist University
in July, 2012. The authors are grateful to two anonymous referees for their comments
and suggestions which greatly improve the presentation of this paper.

\end{document}